\documentclass[10pt]{amsart}

\usepackage{amsthm, amsmath}
\usepackage{amssymb}
\usepackage[all]{xy}
\numberwithin{equation}{subsection}

\theoremstyle{plain}

\newtheorem{theorem}{Theorem}[subsection]

\newtheorem{lemma}[theorem]{Lemma}

\newtheorem{prop}[theorem]{Proposition}

\begin{document}

\title{Corrigendum: The base change fundamental \\ lemma for central elements in \\ parahoric Hecke algebras}

\author{Thomas J. Haines
}

\maketitle

\markright{Corrigendum to \cite{H09}}

\section{Introduction}

In section 2.2 of \cite{H09}, there is a minor misstatement that this note will correct and clarify.  It has no effect on the main results of \cite{H09}, but nevertheless this corrigendum seems necessary in order to avoid potential confusion.  Also, I take this opportunity to point out a related typographical error in \cite{BT2}, section 5.2.4, and to address some matters of a similar nature.

I am very grateful to Brian Smithling and Tasho Kaletha, who informed me that something was amiss in section 2 of \cite{H09}.


\section{Notation}

All notation will be that of \cite{H09}, except for the correction in notation discussed below.

\section{Correction}

In \cite{H09}, section 2.2, the ``ambient'' group scheme $\mathcal G_{{\bf a}_J}$ was incorrectly identified with the group scheme whose group of $\mathcal O_L$-points is the full fixer of the facet ${\bf a}_J$.   In the notation of Bruhat-Tits \cite{BT2}, which I intended to follow in \cite{H09}, the group scheme whose group of $\mathcal O_L$-points is the full fixer of ${\bf a}_J$ is denoted 
$\widehat{\mathcal G}_{{\bf a}_J}$.  The group scheme $\widehat{\mathcal G}_{{\bf a}_J}$ is defined and characterized in this way in \cite{BT2}, 4.6.26-28.  

The group scheme denoted $\mathcal G_{{\bf a}_J}$ is defined in loc.~cit.~4.6.26 (cf. also 4.6.3-6).  In general, it can be a bit smaller than $\widehat{\mathcal G}_{{\bf a}_J}$ (see below).  In \cite{H09}, the symbol $\mathcal G_{{\bf a}_J}$ should be interpreted as this potentially proper subgroup of the full fixer $\widehat{\mathcal G}_{{\bf a}_J}$.

We have, as stated in \cite{H09}, (2.3.2) and (2.3.3), the equalities\footnote{In light of the typographical error in \cite{BT2}, 5.2.4 explained in section \ref{typo}, the reasoning used in \cite{H09} to justify these equalities is correct.}
\begin{align} 
J(L) = \mathcal G^\circ_{{\bf a}_J}(\mathcal O_L) &= T(L)_1 \cdot \mathfrak U_{{\bf a}_J}(\mathcal O_L) \label{2.3.2} \\
\mathcal G_{{\bf a}_J}(\mathcal O_L) &= T(L)_b \cdot \mathfrak U_{{\bf a}_J}(\mathcal O_L). \label{2.3.3}
\end{align}
In general,
$$
\mathcal G^\circ_{{\bf a}_J}(\mathcal O_L) = \widehat{\mathcal G}^\circ_{{\bf a}_J}(\mathcal O_L) \subset \mathcal G_{{\bf a}_J}(\mathcal O_L) \subset \widehat{\mathcal G}_{{\bf a}_J}(\mathcal O_L),
$$
and both inclusions can be strict.  

\section{Clarification of subsequent statements in \cite{H09}}

${\bf 1.}$  Theorem 2.3.1 of \cite{H09} remains valid as stated, but can be slightly augmented: equation (2.3.1) can be replaced by
\begin{equation} \label{thm_2.3.1}
J(L) = {\rm Fix}({\bf a}_J^{\rm ss}) \cap G(L)_1 = \mathcal G_{{\bf a}_J}(\mathcal O_L) \cap G(L)_1 = \widehat{\mathcal G}_{{\bf a}_J}(\mathcal O_L) \cap G(L)_1.
\end{equation}
Cf. \cite{HRa}, Remark 11.

\bigskip

${\bf 2.}$  Contrary to \cite{H09}, line above equation (2.3.2), our ${\mathcal G}_{{\bf a}_J}$ should not now be identified with the scheme $\widehat{\mathcal G}_{{\bf a}^{\rm ss}_J}$ of \cite{BT2}.

\bigskip

${\bf 3.}$  Corollary 2.3.2 of \cite{H09} remains valid, with the same proof. Indeed, when $G_L$ is split we have $T(L)_b = T(\mathcal O_L) = T(L)_1$ and then from (\ref{2.3.2}) and (\ref{2.3.3}) above we see that $\mathcal G^\circ_{{\bf a}_J}(\mathcal O_L) = \mathcal G_{{\bf a}_J}(\mathcal O_L)$. 

\bigskip

${\bf 4.}$  Lemma 2.9.1 of \cite{H09} remains valid as stated, but in the proof (especially in equations (2.9.1) and (2.9.2)) the symbols $\mathcal G_{{\bf a}_J}(\mathcal O_L)$ and $\mathcal G_{{\bf a}^M_J}(\mathcal O_L)$ should be replaced by  $\widehat{\mathcal G}_{{\bf a}_J}(\mathcal O_L)$ and $\widehat{\mathcal G}_{{\bf a}^M_J}(\mathcal O_L)$, respectively.

\section{Example}

It is sometimes but usually not the case that $\mathcal G_{{\bf a}_J}(\mathcal O_L) = \widehat{\mathcal G}_{{\bf a}_J}(\mathcal O_L)$.  The following is perhaps the simplest example where this equality fails\footnote{Brian Smithling and Tasho Kaletha provided me with another example for the split group ${\rm SO}(2n)$.}.  Take $G$ to be the split group ${\rm PSp}(4)$, and let ${\bf a}_J$ denote the non-special vertex in a base alcove.  Then let $\tau$ denote the element in the stabilizer $\Omega \subset \widetilde{W}(L)$ of the base alcove, which interchanges the two special vertices and fixes ${\bf a}_J$.  The element $\tau$ does not belong to the group $\mathcal G^\circ_{{\bf a}_J}(\mathcal O_L) = \mathcal G_{{\bf a}_J}(\mathcal O_L)$ (cf.~${\bf 3}$ above), since $\tau$ does not belong to $G(L)_1$.  On the other hand $\tau \in \widehat{\mathcal G}_{{\bf a}_J}(\mathcal O_L)$ since it fixes ${\bf a}_J$ and $G(L)^1 = G(L)$ (cf. \cite{BT2}, 4.6.28).  

\section{Typographical error in \cite{BT2}, 5.2.4} \label{typo}

Section 5.2.4 of \cite{BT2} contains four displayed equations.  In all of these equations, the ``hats'' should be removed.  The fact that the final displayed equation
$$
\widehat{\mathfrak G}^\natural_{\Omega}(\mathcal O^\natural) = {\mathfrak G}^\circ_{\Omega}(\mathcal O^\natural) \, \mathfrak Z(\mathcal O^\natural)
$$
is incorrect as stated is shown by the Example above (in light of the fact that for a $K^\natural$-split group such as ${\rm PSp}(4)$ the group scheme $\mathfrak Z$ is connected and the right hand side is simply 
${\mathfrak G}^\circ_{\Omega}(\mathcal O^\natural)$).

All of the displayed equations in \cite{BT2}, 5.2.4 become correct when the ``hats'' are removed.

\section{When is $\mathcal G_{{\bf a}_J}(\mathcal O_L) = \widehat{\mathcal G}_{{\bf a}_J}(\mathcal O_L)$?}

Let us assume (for simplicity) that $G$ is split over $L$.  Then the following give two cases where the equality $\mathcal G_{{\bf a}_J}(\mathcal O_L) = \widehat{\mathcal G}_{{\bf a}_J}(\mathcal O_L)$ holds.  Since $G_L$ is split, by Corollary 2.3.2 of \cite{H09} we automatically have $\mathcal G_{{\bf a}_J}(\mathcal O_L) = {\mathcal G}^\circ_{{\bf a}_J}(\mathcal O_L)$.

\begin{lemma} \label{lemma1}
If $G_{\rm der} = G_{\rm sc}$, then $\widehat{\mathcal G}_{{\bf a}_J}(\mathcal O_L) = \mathcal G_{{\bf a}_J}(\mathcal O_L)$.
\end{lemma} 

\begin{proof}
Let $\mathcal I = {\rm Gal}(\overline{L}/L)$ denote the inertia group.  Recall that $G(L)_1$ is the kernel of the Kottwitz homomorphism
$$
G(L) \rightarrow X^*(Z(\widehat{G})^\mathcal I)
$$
and $G(L)^1$ is the kernel of the map 
$$
G(L) \rightarrow X^*(Z(\widehat{G})^\mathcal I)/torsion
$$
derived from the Kottwitz homomorphism.  Our hypotheses imply that $X^*(Z(\widehat{G})^\mathcal I) = X^*(Z(\widehat{G}))$ is torsion-free, and hence $G(L)^1 = G(L)_1$.  But then $\widehat{\mathcal G}_{{\bf a}_J}(\mathcal O_L)$, being by \cite{BT2}, 4.6.28 the fixer of ${\bf a}_J^{\rm ss}$ in $G(L)^1$, obviously coincides with $\mathcal G^\circ_{{\bf a}_J}(\mathcal O_L)$, the fixer of ${\bf a}_J^{\rm ss}$ in $G(L)_1$ (cf. (\ref{thm_2.3.1}) above). 
\end{proof}

\begin{lemma} \label{lemma2}
If the closure of ${\bf a}_J$ contains a special vertex $v$, then $\widehat{\mathcal G}_{{\bf a}_J}(\mathcal O_L) = \mathcal G_{{\bf a}_J}(\mathcal O_L)$.  
\end{lemma}

\begin{proof}
By \cite{BT2}, 4.6.26, we have $\widehat{\mathcal G}_{{\bf a}_J}(\mathcal O_L) = \widehat{N}^1_{{\bf a}_J} \, \mathcal G_{{\bf a}_J}(\mathcal O_L)$, where $\widehat{N}^1_{{\bf a}_J}$ denotes the fixer in $N = N_G(T)(L)$ of ${\bf a}_J$.  Hence, it suffices to show that $\widehat{N}^1_{{\bf a}_J} \subset G(L)_1$.  Let $K = K_v$ be the special maximal parahoric subgroup of $G(L)$ corresponding to $v$, and realize the finite Weyl group $W$ at $v$ as $W = (K \cap N_G(T))/T(\mathcal O_L)$, cf.~\cite{HRa}.  As in loc.~cit.,~the choice of the special vertex $v$ gives us a decomposition of the extended affine Weyl group as $X_*(T) \rtimes W$.  For $n \in \widehat{N}^1_{{\bf a}_J}$ let $t_\lambda w \in X_*(T) \rtimes W$ denote the corresponding element. 

We need to show that $t_\lambda w $ belongs to the affine Weyl group, since such an element will automatically belong to $G(L)_1$, and that would be enough to prove that $n \in G(L)_1$.  We need to show $\lambda$ is in the coroot lattice $Q^\vee$.  But $t_\lambda w$ fixes $v$, that is,
$$
\lambda + w(v) = v.
$$
On the other hand
$$
v - w(v) \in Q^\vee,
$$
since $v$ is a special vertex.  Thus $\lambda \in Q^\vee$ and we are done.
\end{proof}

\section{Comparing Iwahori subgroups over $F$}

The ``naive'' Iwahori subgroup that often appears in the literature (e.g.~\cite{C}, \cite{Mac}), can be identified with the group
$$
\widetilde{I} := G(F) \cap {\rm Fix}({\bf a}^\sigma) = G(F)^1 \cap {\rm Fix}(({\bf a}^{\rm ss})^\sigma).
$$
This contains the group
$$
\widehat{\mathcal G}_{{\bf a}}(\mathcal O_F) = G(F)^1 \cap {\rm Fix}({\bf a}),
$$
(cf.~\cite{BT2}, 4.6.28).  The ``true'' Iwahori subgroup over $F$ is defined to be
$$
I := G(F) \cap (G(L)_1 \cap {\rm Fix}({\bf a}))  = \mathcal G^\circ_{{\bf a}}(\mathcal O_F) 
$$
(see \cite{HRa}) which turns out to have the alternative description
$$
I = G(F)_1 \cap {\rm Fix}({\bf a}^\sigma),
$$
see \cite{HRo}, Remark 8.0.2.
Thus, we always have the inclusions
\begin{equation*}
I \subseteq \widehat{\mathcal G}_{{\bf a}}(\mathcal O_F) \subseteq \widetilde{I}.
\end{equation*}
  
In general, we have $\widetilde{I} \neq I$; for example, in the case of $G = D^\times/F^\times$ we have $\widehat{\mathcal G}_{\bf a}(\mathcal O_F) \neq \widetilde{I}$  (see Remark 8.0.2 of \cite{HRo}).  

\begin{lemma}  Suppose $G$ is split over $L$.  Then $I = \widehat{\mathcal G}_{{\bf a}}(\mathcal O_F)$.
\end{lemma}

\begin{proof}
Use Lemma \ref{lemma2}.
\end{proof}

\begin{prop}
If $G$ is unramified over $F$, then $I = \widehat{\mathcal G}_{{\bf a}}(\mathcal O_F) = \widetilde{I}$.
\end{prop}

\begin{proof}
It is enough to prove $I = \widetilde{I}$.  Let $v_F$ denote a hyperspecial vertex in the closure of $({\bf a}^{\rm ss})^\sigma$, and let $K = K_{v_F}$ denote the corresponding special maximal parahoric subgroup of $G(F)$.  Following \cite{HRo}, define $\widetilde{K} = G(F)^1 \cap {\rm Fix}(v_F)$; recall also that $K = G(F)_1 \cap {\rm Fix}(v_F)$. By loc.~cit., it is clear that when $G$ is unramified over $F$ we have $\widetilde{K} = K$.   
On the other hand, the inclusion
$\widetilde{I} \subset \widetilde{K}$ clearly induces an injection
$$
\widetilde{I}/I \hookrightarrow \widetilde{K}/K.
$$
Thus $\widetilde{I}/I$ is trivial.
\end{proof}

\small
\bigskip
\obeylines
\noindent
University of Maryland
Department of Mathematics
College Park, MD 20742-4015 U.S.A.
email: tjh@math.umd.edu

\end{document}